\documentclass[12pt,leqno]{amsart}
\usepackage{amsmath}
\usepackage{amssymb}
\usepackage{amsthm} 
\usepackage{epsf}
\usepackage{graphicx}
\usepackage{esint} 
\usepackage{dsfont}
\usepackage[pagebackref]{hyperref} 
\hypersetup{backref,  pdfpagemode=FullScreen,  colorlinks=true} 
\usepackage{color}
\usepackage{enumerate}

\numberwithin{equation}{section}

\newcommand{\ms}{\medskip}

\newcommand{\R}{\mathbb{R}}
\renewcommand{\H}{\mathcal H}

\newcommand{\bC}{\mathbb {C}}

\renewcommand{\d}{\partial}
\newcommand{\dist}{\,\mathrm{dist}}

\newcommand{\sm}{\setminus}

\newcommand{\wt}{\widetilde}

\newcommand{\ol}{\overline}
\newcommand{\ub}{\underbar}

\newcommand{\cQ}{{\mathcal Q}}

\newcommand{\fS}{{\mathfrak S}}

\theoremstyle{plain}
\newtheorem{theorem}[equation]{Theorem}

\theoremstyle{definition}

\newtheorem{example}[equation]{Example}

\theoremstyle{remark}
\newtheorem{remark}[equation]{Remark}

\newcommand{\dv}{\operatorname{div}}

\newcommand{\RR}{{\mathbb{R}}}

\definecolor{greenish}{cmyk}{0.91,0,0.88,0.12}

\newcounter{hcountno}

\begin{document}
\title{Good elliptic operators on Cantor sets}

\author[David]{Guy David}
\address{Guy David. Universit\'e Paris-Saclay, CNRS, Laboratoire de math\'ematiques d'Orsay, 91405 Orsay, France}
\email{guy.david@universite-paris-saclay.fr}

\author[Mayboroda]{Svitlana Mayboroda}
\address{Svitlana Mayboroda. School of Mathematics, University of Minnesota, Minneapolis, MN 55455, USA}
\email{svitlana@math.umn.edu}

\thanks{David was partially supported by the European Community H2020 grant GHAIA 777822, 
and the Simons Foundation grant 601941, GD. Mayboroda was supported in part by the Alfred P. Sloan Fellowship, 
the NSF grants DMS 1344235, DMS 1839077, and the Simons foundation grant 563916, SM}
\date{}

\maketitle

\begin{abstract}
It is well known that a purely unrectifiable set cannot support a harmonic measure 
which is absolutely continuous with respect to the Hausdorff measure of this set. 
We show that nonetheless there exist elliptic operators on (purely unrectifiable) 
Cantor sets in $\R^2$ whose elliptic measure is absolutely continuous, and in fact,  
essentially  proportional to the Hausdorff measure.
\end{abstract}

\ms\noindent{\bf Key words:}  Harmonic measure, Green function, Elliptic operator,
Counterexample, Cantor set.

\ms\noindent
AMS classification 35J15, 35J08, 31A15, 35J25.

\tableofcontents

\section{Introduction}
Starting with the seminal 1916 work of F. and M. Riesz \cite{RR}, 
considerable efforts over the century culminating in the past 10-20 years, 
identified necessary and sufficient geometric conditions on the domains for which harmonic measure is absolutely 
continuous with respect to the Hausdorff measure of the boundary. 
In particular, it was established that purely unrectifiable sets are those for which 
harmonic measure is necessarily singular with respect to the Hausdorff measure of the boundary of the domain. 
Moreover, it was shown that a similar statement holds for operators reasonably close to the Laplacian. 

It turns out, however, that even for a purely unrectifiable set there may exist  
a ``good" elliptic operator. 
The main result of this paper is the construction of elliptic operators $L$
associated to the four corners Cantor set $K$ of dimension $1$ in the plane,
whose elliptic measure $\omega_L$ is essentially proportional to the 
one-dimensional Hausdorff measure on $K$. 
To the best of our knowledge, this is the first result of this nature.

We shall concentrate mainly on the emblematic Garnett-Ivanov
Cantor set of dimension~$1$ in the plane, also known as the four-corners Cantor set
(see Figure \ref{f2}),  
because it is probably the most celebrated example of a one-dimensional, 
Ahlfors regular set, such that the harmonic measure and the Hausdorff measure 
$\H^1$ are mutually singular on $K$, but we will explain in Section \ref{S5}
that our construction is fairly flexible.

The constructed operators $L=\dv A \nabla$ are of divergence form. 
Moreover, they are scalar, that is, we can write them as 
\begin{equation} \label{1.1}
L = \dv a \nabla,
\end{equation}
where $a$ is a 
continuous
scalar function on $\R^2$ (as opposed to a general matrix-valued function $A$)
such that 
$C^{-1} \leq a(X) \leq C$ for $X\in \Omega = \R^2 \sm K$. 

Even though for the Cantor set the question whether there exists an elliptic operator for which elliptic measure is absolutely continuous with respect 
to the Hausdorff measure on the boundary was also open for general elliptic 
operators $L = \dv A \nabla$, we aimed to have a solution in the smaller class
of isotropic operators $L = \dv a \nabla$ as above
which seem to be more clearly geometrically relevant. Unfortunately, this is also harder: 
if one agrees to work with a general case of matrix-valued $A$ the theory of quasiconformal mappings 
becomes an ally.

Let us briefly discuss the history of known results pertaining to the following general question. 
Given a domain $\Omega \subset \R^n$, bounded by $E = \d\Omega$,  and a divergence form operator 
$L =  \dv A \nabla$, where 
$A= A(X)$ is an elliptic matrix defined for $X \in \Omega$, when can we say that
the elliptic measure $\omega_L$ associated to $E$ (on $\Omega$) is absolutely
continuous with respect to the surface measure on $E$ or more generally, when
$E$ is not assumed to be smooth,  
to the relevant 
Hausdorff measure or an Ahlfors regular measure
$\sigma$ living on $E$? 

Consider first the case of the Laplacian. In the positive direction,
many results give the absolute continuity of $\omega_\Delta$ when $E$ is rectifiable and 
some sort of topological connectedness condition is satisfied;
see for instance \cite{RR, L, KL} in the plane, and \cite{Dah, DJ, Se, HM1} 
in higher dimensions,
and for instance \cite{Bad} in the absence of Ahlfors regularity. In the converse direction, 
it was long known that the harmonic measure on the Garnett-Ivanov Cantor $K$
is singular with respect to the natural measure, but the precise more general
results are much more recent. In particular, rectifiability was identified as a necessary condition for the absolute continuity of harmonic measure only in 2016  \cite{AHM3TV}, but even then the exact necessary topological (connectedness) assumptions remained elusive. Finally, in 2018, the sharp necessary and sufficient conditions were established in \cite{AHMMT}. Many of these results generalize to operators $L$ other than the Laplacian 
but morally close to it, for instance, those satisfying a suitable square Carleson measure condition.
See \cite{KePiDrift, HMMTZ}.

Concerning Cantor sets, it has been known since \cite{Car} that 
in the case of $K$, the (usual) harmonic measure $\omega_\Delta$ is singular 
with respect to the Hausdorf measure, and even that it lives on a subset of dimension
strictly smaller than $1$. The same thing is true even 
in larger dimensions, and Carleson's result was later generalized to larger classes of fractal
sets; see for instance \cite{Bat, BZ}. These results are 
quite delicate (and the condition that $\omega_\Delta$ is supported by a set of dimension $<1$ is also 
significantly stronger than the  
mere singularity). See also \cite{Azzam} for a more recent result with much less structure.

The present papers deals with operators $L$ that are {\it{not}} close to the Laplacian.
For its authors, the issue arose when they tried to define reasonable elliptic operators
on $\R^n \sm E$, where $E$ is an Ahlfors regular set of dimension $d < n-1$. In such a 
context, it was shown in \cite{DFM1, DFM2} that $L$ should rather be of the form $L = \dv A \nabla$,
where the product of $A$ by $\dist(X,E)^{n-d-1}$ satisfies the usual ellipticity conditions.
These operators cannot be close to the Laplacian (except in spirit). 
What was perhaps more surprising, is that in order to prove the absolute continuity of the  
elliptic measure, we had to work with a very particular choice of coefficients, not the one driven by the powers 
of the Euclidean distance. And indeed, if 
$$
D_\alpha(X)=\left(\int_E |X-y|^{-d-\alpha} \,d\H^d_{\vert E}(y)\right)^{-\frac 1\alpha}, 
\quad X\in \RR^n\setminus E, \quad \alpha>0,
$$
and $L_\alpha=\dv D_\alpha(X)^{-n+d+1} \nabla$ then we can prove 
that the elliptic measure is absolutely continuous (and given by an $A_\infty$ weight)
when $E$ is uniformly rectifiable; see \cite{DFM3, DM, Fen}.

For these results the converse is not known. In fact, the authors of \cite{DEM} have discovered a ``magical" counterexample which really brings us close to the subject of the present paper. 
It turns out that when $d<n-2$ and $\alpha =n-d-2$, the elliptic measure for the operator $L_\alpha$ defined above is absolutely continuous with respect to the Hausdorff measure for {\it any Ahlfors regular set of dimension $d$.} 
The dimension does not have to be an  
integer, and the set in question does not have to be rectifiable or carry any other geometric characteristics on top of Ahlfors regularity. Of course, in this case the coefficients depend on the  set via a particular  distance function.

In view of all these results the following question is quite natural. 
Even in the classical case of sets with $n-1$ dimensional boundary, 
given a bad (totally unrectifiable set) like $K$ above, is it possible to find elliptic operators $L = \dv A \nabla$ 
on $\R^n \sm K$ such that $\omega_L$ is nonetheless absolutely continuous (possibly even proportional) 
to the natural measure on $K$?

The present paper shows that the answer is yes for the Cantor set $K \subset \R^2$ described near \eqref{3.1},
and quite a few other ones. In fact, we establish the 
stronger result that the elliptic measure and Hausdorff measure are roughly proportional, in the sense that
if we take a pole $X \in \R^2 \sm B(0,1)$ (that is, far enough from $K$), then
\begin{equation} \label{1.2}
C^{-1} \H^1_{\vert K} \leq \omega^X \leq C  \H^1_{\vert K}.
\end{equation}

See Theorem \ref{t41} for the precise statement, and Section \ref{S5} for 
other sets $K$.

The basic idea is very simple: we shall be able to construct the Green function
with a pole at $\infty$. Usually, one does not dream of computing the Green function explicitly,
except in the very simple cases of the Laplacian on a disc or a half-space, 
or when some unexpected miracle happens (see, e.g., the aforementioned example \cite{DEM}). 
But here the situation is different: we construct the Green function $G$ first, 
subject to suitable constraints, and then compute the coefficients of $L$
in terms of $G$. This will require some care because we want $G$ to be a solution of 
$\dv a \nabla G = 0$ with $C^{-1} \leq a \leq C$, 
and $a$ is very far from unique given $G$,
but we have a fair chance.

We shall use the fact that we work in dimension $2$ because $G$, and locally a
conjugate function, will be computed from their level lines and the way they cross.
In particular, we will start from a family of level lines and construct the orthogonal curves,
and this is certainly easier in the plane. Also, we are not sure that there
is a good enough notion of conjugate function in higher dimensions. 
Of course we can create examples in $\R^n$, for instance on the product of $K$ with $\R^{n-2}$,
but in terms of construction this is obviously cheating.

In the course of the proof we discovered a number of the properties of the equation $Lu = 0$,
with $L$ as in \eqref{1.1}, in dimension $2$, concentrating on a somewhat less conventional direction from solution to the coefficients of the PDE. Some of these are probably very well known, 
particularly in connection with the so-called Calder\'on inverse problem, but 
they were of considerable help to us for understanding how to find coefficients $a$
such that $\dv a \nabla G = 0$ for a given function $G$. We explain this in Section~\ref{S2}. 

The definition of $K$ and main construction, with the level curves and the pictures, 
will be done in Section~\ref{S3}. 
As we just said, the level sets of $G$ will be constructed so that 
$G$ is a solution of some equation $LG=0$, but we will need to make sure that the function
$a$ computed from $G$ has uniform bounds. For this, the fractal nature of $K$ will
be useful, because it will allow us to prove the desired estimates only on some
fundamental domains, and then glue the different pieces.

We shall then explain, in a short Section \ref{S4} the relation between the 
constructed Green function and the elliptic measure associated to $L$; 
this is also the section where our main theorem will be stated (hopefully with no surprise).

In the last Section \ref{S5}, we explain that our construction is not so rigid. For instance it works
essentially with no modification for a non-fractal variant of $K$ where we are also allowed to rotate the
squares independently. 
This is interesting, because for these Cantor sets, as far as we know the fact that
the harmonic measure $\omega_\Delta$ 
lives on a set of smaller dimension 
has been established only recently, with a very complicated proof, while the case of 
$K$ was treated quite some time ago \cite{Car}, 
using the fractal nature of $K$. 

For the operator $L$ that we construct, the Green function $G(X)$ is equivalent
(for $X$ close to $K$ and the other pole far away) to the distance $\dist(X,K)$; 
it is amusing that for slightly different functions $a$ (but this time degenerate elliptic), 
or, alternatively, for similar Cantor sets of different dimensions, we get Green functions $G$
that are equivalent to different powers of $\dist(X,K)$ 
(and $\omega_L$ is still proportional to the Hausdorff measure).
This is related to the invariance properties of the equation
$Lu = 0$, with $L$ as in \eqref{1.1}, but also to the reason why we
prefer to have a scalar function $a$.

As we mentioned above, another reason for emphasizing scalar coefficients (and for the difficulties this entails) 
comes from the theory of quasiconformal mappings.  
Consider the similar question where $K$ is instead a snowflake of dimension $\alpha > 1$
in $\R^2$. If you want to find an elliptic operator $L = \dv A \nabla$ such that
the associated elliptic measure $\omega_L$ is absolutely continuous
with respect $\H^{\alpha}_{\vert K}$, 
this is in fact (too) easy: we use a
quasiconformal mapping $\psi : \R^2 \to \R^2$ that maps 
the line $\ell$ (or a circle, depending on whether $K$ is unbounded or bounded)
to $K$, and then use $\psi$ to move the Laplacian on a component of $\R^2 \sm \ell$
to the desired component of $\R^2 \sm K$. It is known that the
conjugated operator is an elliptic
operator $L = \dv A \nabla$, and the absolute continuity of $\omega_L$
with respect to $\H^\alpha$ follows directly from the corresponding absolute continuity
result for $\ell$ and $\Delta$.

Of course in the case of the Cantor set $K$, we cannot do that, 
even if we allow general elliptic matrices $A$,
but this says that the class of elliptic operators is 
really too large. As far as the authors know, the same question for a snowflake,
but with a scalar operator $L$ as in \eqref{1.1}, is open. Possibly the ideas in this paper 
could help, but one would need to write down level sets in a very careful way.

\section{About the equation $\dv a \nabla u = 0$ in the plane}
\label{S2}

In this section we try to see how a given function $u$ on the plane
can be seen as a solution of an equation $Lu = 0$, with $L$ as in \eqref{1.1}
and a function $a$ that we could compute in terms of $u$.

We want to do this in a geometrical way, in terms of the level sets of $u$,
by introducing a second function $v$, which is related to $u$ (as is conventional, 
we will say conjugated in analogy with harmonic functions)
and satisfies a similar equation but with the function $a^{-1}$. 
We proceed locally in an open set $U$, where we assume that $u$ is of class $C^4$ 
and with $\nabla u \neq 0$. 
There could be a brutal analytic way to find $v$ and $a$, but the geometry 
of level sets seems much easier to understand.

Associated to $u$ are the level sets $\gamma_s = u^{-1}(\{ s \})$,
and its nonvanishing gradient. We can use the gradient as a vector field, and then
solve the equation $z'(t) = \nabla u(z(t))$ to get a family of orthogonal curves
$\Gamma_\theta$. To be more precise, we first parameterize
one of the level curves, say, $\gamma_0$, call $\theta \to \gamma_0(\theta)$
this paremeterization, and then solve the equation 
$z'(t) = \nabla u(z(t))$ with the initial condition $z(0) = \gamma_0(\theta)$
to get the curve $\Gamma_\theta$. We can extend $\Gamma_\theta$
in both direction, as long as it stays in the domain $U$,
and we know from the uniqueness of solutions that the curves $\Gamma_\theta$
never cross. They may be periodic, though.

Locally they cover the space, in the sense that if
$z$ lies in some $\Gamma_{\theta_0}$, it is easy to see (by running along the vector field
backwards) that every point of a small ball around $z$ lies in the union of the 
$\Gamma_\theta$, $\theta$ close to $\theta_0$.
Indeed, we can run along $\Gamma_{\theta_0}$ backwards starting from
$z$ (and up to the point $\gamma_0(\theta_0)$.
For $z'$ close to $z$, we can run the same vector field, and get a 
curve that stays close to $\Gamma_{\theta_0}$. In fact, the solution is a $C^2$ 
function of $z$, because $\nabla u$ is of class $C^3$ 
(we are not aiming for optimal regularity here). Also, since 
$\Gamma_{\theta_0}$ meets $\gamma_0$ transversally at 
$\gamma_0(\theta_0)$, we can apply the implicit function theorem to 
prove that there is a unique $\theta$, $\theta$ close to $\theta_0$,
that contains $z'$. Thus we get a $C^2$ mapping $v$, defined near $z$, which to $z$ associates 
the unique $\theta$ near $\theta_0$ such that $z' \in \Gamma_{\theta}$. In fact, the uniqueness is global
(as long as the parameterization of $\gamma_0$ is injective; we would need
to be more specific when $\gamma_0$ is a loop), because 
the different $\Gamma_\theta$ do not cross. Finally, 
$\nabla v(z) \neq 0$, because it is obtained from the implicit function theorem.

At this point, we have, in an open set that contains $\gamma_0$, two functions $u$ and
$v$ that can be used as coordinates. The pair $(u,v)$ satisfies the desirable orthogonality 
relation
\begin{equation} \label{2.1}
\nabla u(z) \perp \nabla v(z).
\end{equation}
This does not mean that the change of variables defined by $(u,v)$ is conformal,
because although the gradients (or the level lines) are orthogonal, the two lengths
$|\nabla u(z)|$ and $|\nabla v(z)|$ are different. In fact, if we wanted to be sure that
$(u,v)$ defines a conformal change of variables, we should require $u$ to be harmonic
(and then we'll see that $v$ is harmonic too). So for example,
replacing $(u,v)$ with $(u+v, u-v)$ gives another pair that does not satisfy 
\eqref{2.1} in general. 

Our next step is to use this relation to show that,
in any domain where $u, v$ are functions of class $C^2$ in an open set, where 
$\nabla u \neq 0$ and $\nabla v \neq 0$, and \eqref{2.1} holds,
\begin{equation} \label{2.2}
\dv a \nabla u = 0, \text{ with } a(z) = \frac{|\nabla v(z)|}{|\nabla u(z)|}.
\end{equation}
Of course the same computation will show that $v$ satisfies
the twin equation $\dv \frac1a \nabla v = 0$.

Let us write this in coordinates and check this near a point $z$,
but avoid to mention the argument $z$ when we don't need to.
The vector $w = (\frac{\d v}{\d y}, -\frac{\d v}{\d x})$
is perpendicular to $\nabla v$ and has the same length as $\nabla v$,
so since $\nabla u$ is proportional to $w$ (by \eqref{2.1}),  
in fact $a\nabla u  = \pm a w$ (because the two vectors have the same length).
The sign is locally constant, so we may now compute
\begin{equation} \label{2.3}
\dv a \nabla u = \dv w = \frac{\d^2 v}{\d x\d y} - \frac{\d^2 v}{\d x \d y}=0,
\end{equation}
as promised.

Here we shall be interested in making sure that
$C^{-1} \leq a \leq C$, or in other words that 
$C^{-1} |\nabla u| \leq |\nabla v| \leq C |\nabla u|$. 
In terms of level curves for $u$ and $v$, assuming that a nice parameterization
of $\gamma_0$ was chosen, the gradient of $G$ is (locally) proportional
to the inverse of the distance between the level sets (divided by an increment 
$\delta_s$), and similarly the gradient of $v$ corresponds to the inverse of the distance between the 
$\Gamma_\theta$ 
(divided by a $\delta_\theta$); 
then we should compare these quantities and say whether they stay between constants. On the picture, it means that if we use roughly equal increments $\delta_x$ and $\delta_\theta$, we should see small rectangles that are not too thin.
I particular, it is all right if the two types of level sets become too sparsely or too densely
spaced, provided that they do this in essentially the same way.

\begin{example} \label{e24}
Let us even describe an example where $u$ and $v$ are conjugated harmonic functions, so that
we can even take $a=1$. Identify $R^2$ with $\bC$, write $z = x+i y$,
also use polar coordinates, and take
\begin{equation} \label{2.5}
u(x,y) =  {\mathcal Re}(z^4) = r^4 \cos 4\theta
\text{ and } 
v(x,y) = {\mathcal Im}(z^4) = r^4 \sin 4\theta.
\end{equation}
It is known (and easy to compute) that $u$ and $v$ are harmonic and
$|\nabla u(x,y)| = |\nabla u(x,y)| = 2|z|^2$.
Of course we are not exactly in the case described above, because $u$ and $v$
have a critical point at the origin, but we do not have to divide by $|\nabla u|$
to know that $u$ and $v$ are harmonic. See Figure \ref{f1} 
for a clumsy attempt to describe the level sets of $u$ (in green) and its gradient lines
(in red). Notice that at the origin, $u$ increases at maximal speed
when $\cos 4\theta = 1$, i.e., when $\theta = k \pi/2$, along the axes, and
$u$ decreases at maximal speed when $\cos 4\theta = -1$, along the diagonals.      
All these lines are special red curves (oriented differently) and are separated
by green lines where $4\theta = \frac{\pi}{2}+k\pi$.

\begin{figure}
\centering
\includegraphics[width=5.cm]{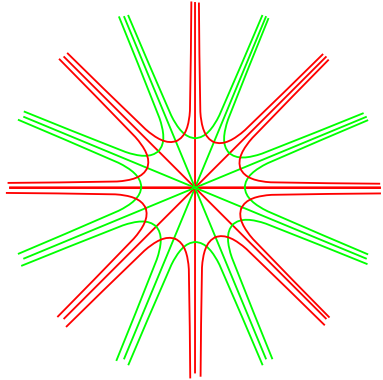}
\caption{Level lines for $u(r,\theta) = r^4 \cos 4\theta$ (in green), and 
the conjugate function $v$ (in red)}
\label{f1}
\end{figure}

The advantage of this example is that the uniform bounds on $|\nabla u|^{-1} |\nabla v|$
are obvious. We may always replace $u$ and $v$ with the new functions $G$ and $R$, where
\begin{equation} \label{2.6}
G(z) = G(0) + \beta u(z) \ \text{ and } R(z) = R(0) + \beta v(z),
\end{equation}
where we can choose $G(0)$, $R(0)$, and $\beta > 0$ as we like
(we prefer $\beta > 0$ because this way we preserve the 
direction of the arrows). This does not change the level lines, just the way they are 
labelled. 
\end{example}

\begin{remark} \label{r27}
If we are given the level lines of $u$, we can deduce the direction of the
gradient, so we can also draw the level lines for $v$. We do not change the picture when we relabel
the level lines of $u$, i.e., compose with a function $f$ and replace $u$ with $w=f\circ u$,
or when we change the parameterization of $\gamma_0$, i.e., in effect, replace
$v$ by some $g \circ v$.

Another way to see this is to observe that when $u$ satisfies $\dv a \nabla u = 0$, then
$w=f\circ u$ satisfies the related equation $\dv b \nabla w = 0$, with $b(z) = a(z)/f'(u(z))$, 
just because $b(z) \nabla(f\circ u)(z) = b(z) f'(u(z)) \nabla u(z) = a(z)\nabla u(z)$,
and then the divergence is the same. 

This is an interesting flexibility that we have with the equation 
$\dv a \nabla u = 0$, at least if we allow ourselves to play with $a$.

We could also replace the function $v$ with another function $g\circ v$; 
this does not change the level sets of $v$, nor the orthogonality condition \eqref{2.1}, 
but as before it changes $a = \frac{|\nabla v|}{|\nabla u|}$.

This last remark is important because it helps us understand that given $u$, we have a large 
choice of functions $a$ such that $\dv a \nabla u = 0$, but some of them are equivalent for
geometric reasons. Even if $u$ is harmonic, we can find lots of operators $L = \dv a \nabla$
such that $Lu=0$, in particular $L = \Delta$ and $L = 2\Delta$, but not only.
The geometry will help us choose good functions $a$.
\end{remark}

\section{Level sets of $G$}
\label{S3}

\begin{figure} [!h] 
\centering
\includegraphics[width=4.cm]{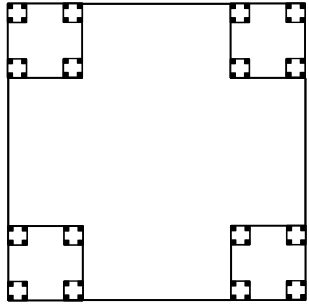}
\caption{The set $K_3$ (three generations of the construction of $K$)}
 \label{f2}
\end{figure}

Let us first define our Cantor set. 
We we start with the square $K_0 = [-\frac12,\frac12] \times[-\frac12,\frac12]$ 
of sidelength $1$. We replace $K_0$ with a set $K_1$, composed of four squares 
$K_{1,j}$ contained in $K_0$, situated at the four corners, and of sidelength $4^{-1}$.
Then we define the sets $K_k$, $k \geq 1$ by induction, to be the set 
$K_k \subset K_{k-1}$, composed of $4^k$ squares of sidelength
$4^{-k}$, and obtained by replacing each of the squares of sidelength 
$4^{1-k}$ that compose $K_{k-1}$ by four squares situated at the four
corners. See Figure \ref{f2}. Our final set is
\begin{equation} \label{3.1}
K = \cup_{k \geq 0} K_k \subset K_0.
\end{equation}
It is a compact set of dimension $1$, such that $0 < \H^1(K) < +\infty$,
which is also Ahlfors regular and has been known in particular to be a simple example
of a compact set such that $\H^1(K) > 0$ and with a vanishing analytic capacity;
see \cite{Ga, Iv}.  It is also known that $\H^1_{\vert K}$ and the harmonic measure on $K$ 
are mutually singular. 

Our goal for this section is to construct a function $G$ on $B_0 = B(0,1)$, 
which we will decide is the restriction to $B_0$ of the desired Green function. 
We can already decide that
\begin{equation} \label{3.2}
\text{$G(z) = 0$ on $K$, and $G(z) = 1$ on $\d B(0,1)$.}
\end{equation}

We will also define $G$ in a self-similar way, so we consider the four centers
$z_j$ of the four squares $Q_j$ that compose $K_1$. Take for instance
$z_1 = (\frac{7}{16},\frac{7}{16})$ and turn in the trigonometric direction, i.e.,
take $z_2 = (-\frac{7}{16},\frac{7}{16})$, $z_3 = (-\frac{7}{16},-\frac{7}{16})$,
and $z_4 = (\frac{7}{16},-\frac{7}{16})$. Then, for $i = 1, \ldots 4$,
set $B_j = B(z_j,1/4)$. Thus $B_j$ is the analogue for $Q_j$ of 
$B_0$ for $K_0$. The radius $r_0=1$ was chosen sufficiently large for $K_0$
to be contained in $B_0$ (and $Q_j$ in $B_j$), but small enough for
$B_j$ to stay away from the axes, i.e., to be contained in the same quadrant as $Q_j$.
We will in fact construct $G$ on the annular region
\begin{equation} \label{3.3}
A_0 = B_0 \sm \Big[\bigcup_{j=1}^4 B_j \Big],
\end{equation}
and in order to be able to glue easily, we will make sure that
\begin{equation} \label{3.4}
G(z) = \frac14 \text{ for } z \in \bigcup_{j=1}^4 \d B_j 
\end{equation}

\begin{figure} [!h]  
\centering
\includegraphics[width=4.cm]{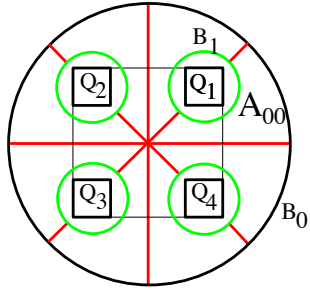}
\caption{The cubes $Q_j$, the balls $B_0$ and $B_1$, and the annuli $A_0$ 
(in the big ball, outside of the green ones) and $A_{00}$ (one eighth of the picture)
}
\label{f3}
\end{figure}

In addition, let us use the symmetries of our set. Call $\fS$ the collection
of symmetries with respect to the two axes and the two diagonals;  
we decide that
\begin{equation} \label{3.5}
G(\sigma(z)) = G(z) \text{ for $z\in A_0$ and }\sigma \in \fS;
\end{equation}
and because of this it will be enough to define $G$ on the smaller region 
\begin{equation} \label{3.6}
A_{00} = \big\{ (x,y) \in A_0 \, ; \, x \geq 0, y \geq 0, \text{ and } x \geq y \big\}.
\end{equation}
See Figure \ref{f3} 
for a first sketch of the level lines of $G$, in green, and its gradient lines 
(or the level lines of the conjugated function), in red. 
Observe before we start that because of the symmetries, we need to have a critical point
at the origin; the part of the axes and the diagonals that lie in $A_{0}$
are really gradient lines of $G$, and they part at the critical point $0$.

We will fill the picture little by little; we will need to pay a special attention to what happens 
near the critical point, because this is the place where it is not obvious that our curves have a
uniform behavior.

We start with the definition of $G$ near the boundary circles, where we slightly prefer
$G$ to be harmonic. So pick a radius $\rho$ just a little bit smaller than $r_0=1$, and set
\begin{equation} \label{3.7}
G(z) = 1 + \ln(|z|) \text{ for $z \in B_0 \sm B(0,\rho)$.}
\end{equation}
This respects the symmetries and \eqref{3.2}; in principle, we should only have taken
\eqref{3.7} in $A_{00} \cap B_0 \sm B(0,\rho) = A_{00} \sm B(0,\rho)$, but this is the same.
We do the same thing near the $B_j$, i.e., choose a radius $\rho_1$ just a tiny bit
larger than $1/4$ and set
\begin{equation} \label{3.8}
G(z) = \frac14 + \frac14\ln(4|z-z_1|) \text{ for $z \in A_{00} \sm B(z_1,\rho_1)$}
\end{equation}
the other balls $B_j$ would be taken of by symmetry. 

\begin{figure} [!h] 
\centering
\includegraphics[width=5.cm]{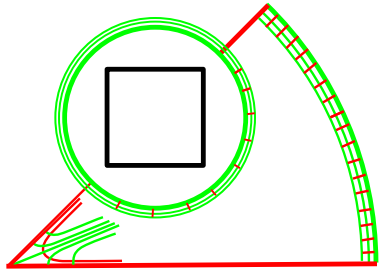}
\caption{The function $G$ (level lines in green, gradient lines in red) in the 
three special regions. The lower left part is a copy of part of Figure \ref{f1}}
 \label{f4}
\end{figure}

Next we take care of the situation near the origin. We decide that in a small ball
$B(0,\rho_2)$, we take the functions $G$ and its conjugate $R$ according to the
formula \eqref{2.6}, where we can still choose $G(0) \in (\frac14,1)$ and
$R(0)$ so that the picture looks nicer (we shall return to this issue when we glue).
In the meantime we can at least draw the green curves and their orthogonal red curves
in $B(0,\rho_2)$ (without labelling them yet). This gives something like Figure \ref{f4}
when we restrict to $A_{00}$ (as before our choice of $G$ preserves the symmetry). 
In particular, \eqref{3.4} is satisfied.

So we have already drawn the level sets of $G$ and its conjugated function $R$
(i.e., the gradient lines) in three small regions of $A_{00}$; 
now we complete the red lines, subject to the following constraints, that we claim are easy
to implement:
\begin{equation} \label{3.9}
\text{ the unit vector field $V_R$ tangent to the red lines is smooth on $A_{00} \sm B(0,\rho_2/2)$}
\end{equation}
(there is indeed a singularity at the origin, but it is well controlled), and in addition
\begin{equation} \label{3.10}
\text{ the intersection of the diagonal and the $x$-axis with $\ol A_{00}$ are red curves,}
\end{equation}
and even slightly more, the red vector field $V_R$ is parallel to the straight parts of the 
boundary of $A_{00}$ in the neighborhood of the corresponding parts of the boundary, and consequently 
when we derive the Green lines, they will meet the straight part of the boundary transversally, and of course 
(since we mentioned vector fields), the red lines do not cross. 
Finally, we need to 
patch the different red lines, for instance in the zone below
and to the right of $B_1$, in such a way that
\begin{equation} \label{3.11}
\begin{aligned}
&\text{the red curve $\Gamma_\theta$ that leaves from $\d B_1$ at the
point $z_1 + \frac14 e^{-i\frac{\pi}{4}+i\theta}$, $0 \leq \theta \leq \pi$,}
\cr&\, \hskip 6cm
\text{leaves $B_0$  at the point $e^{i\theta/4}$.}
\end{aligned}
\end{equation}
That is, we have a mapping from a part (a half) of $\d B_1$ to a part of $\d B_0$ (one eighth), 
obtained by associating the two endpoints of each red curve, and we want it to be not only bijective, 
but run at constant speed. Here the two arcs have the same length, so this also means that the mapping is 
locally isometric (for the arclength).
This is a little unpleasant to do practically, 
so let us reduce this to a simpler problem. We use a mapping $\psi$ to send $A_{00}$
to a rectangle $R$, so that the two circular pieces of the boundary,
namely $\d B_1 \cap \d A_{00}$ and $\d B_0 \cap \d A_{00}$, are sent to
two opposite edges of $R$, which we call $I_1$ and $I_3$, and even at constant speed $1$.
The two arcs have the same length, so this is possible. We can also make sure
that the mapping is a smooth diffeomorphism (with a bounded derivative for the inverse),
except at the origin where we can make it equal to $z \mapsto z^4$. 

We already constructed (pieces of) our red curves $\Gamma_\theta$ in a small neighborhood 
in $\ol A_{00}$ of $\d A_{00}$, and we consider their images $\wt\Gamma_\theta
= \psi(\Gamma_\theta)$ in $R$. First of all, notice that the singularity at $0$
of our collection of red curves disappears, so we really have pieces of curves that are smooth
(where they are defined) and leave $I_1$ and arrive on $I_3$ perpendicularly.
Those that are globally defined, along the remaining edges $I_2$ and $I_3$, also do not cross, 
and follow nicely the edges. In addition, we may modify them a little near the
middle of $I_2$ and $I_4$, if needed, so that the the endpoint $y(\theta)$ of 
$\wt\Gamma_\theta$, which is defined for $\theta$ near $0$ and near $\pi$, goes along 
$I_3$ at constant speed. Along $I_1$ and $I_3$, our red curves touch the boundary
perpendicularly. By restricting to a smaller neighborhood of $\d R$
(or equivalently $\d A_{00}$), we can make sure that they are graphs of Lipschitz functions with a small constant (assuming that $I_1$ and $I_3$ are vertical). 
At this point we claim that it is easy to extend the curves so that they fiber
$R$ and $y(\theta)$ runs at constant speed $1$ between the two vertices of $I_3$.
Then the inverse images by $\psi$ solve our initial problem.
We kindly leave the verification to the reader.

\begin{figure} [!h]  
\centering
\includegraphics[width=5.cm]{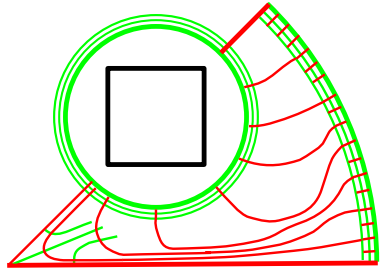}
\caption{The completed red lines, subject to the constraints above. We make them turn
in the region in the right so that the two arcs of circle connect at equal speeds}
\label{f5}
\end{figure}

\begin{figure} [!h] 
\centering
\includegraphics[width=5.cm]{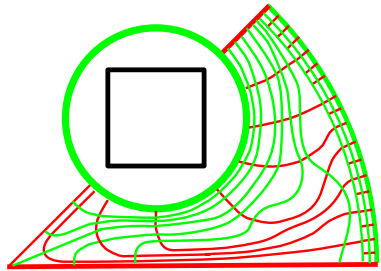}
\caption{The level lines of $G$ (in green). 
Notice that a principle once stated by B. Dahlberg verifies: in spite of the existence and regularity theorem, it is very hard to draw orthogonal green curves so that they look nice.} 
 \label{f6}
\end{figure}

\ms
When all this is done, we get a picture like Figure \ref{f5}. 
The collection of red curves, labeled by $\theta$, defines a function $R$
on the region $A_{00}$. The green curves that are drawn in  Figure~\ref{f6}, 
are just the orthogonal curves, obtained as in Section \ref{S2} by following a vector field
$V_G$ perpendicular to the previous one. In  
$B(0,1/4)$ we already had what we needed,
so we don't need to solve a singular vector field problem, just notice that the two definitions patch.
The smoothness of $V_R$ gives the existence and regularity of the green lines. A continuity argument shows that they start from the upper part $D_1$ of the first diagonal, perpendicularly
to $D_1$, and then end up along the rest of the straight boundary, on the union of the
lower part $D_2$ of the 
first diagonal, followed by the piece of the $x$-axis. 
The perpendicular landing comes from \eqref{3.10}, the rest comes from the existence, uniqueness, and smoothness for the integral curves. 

\begin{remark} \label{r31}
Finally the reader may wonder why we asked for a matching condition on the extremities 
of the red curves. Near $\d B_0$ and $\d B_1$, we decided a formula for $G$,
and our construction ensures that $R$ is a conjugate function to $G$, as in the orthogonality
condition \eqref{2.1}. But there is a dangerous pitfall here, which we try to explain so that
the reader does not make the same mistake as some of the authors.

It is true that $G$ is harmonic near $\d B_0$ and $\d B_1$, so we know for sure
that it satisfies $LG=0$, with $L=\Delta$ and $a=1$. But many other choices of $a$
are possible, for instance $a=1$ near $\d B_0$ and $a=12$ near $\d B_1$, not to mention
more exotic choices. Here it is important that $a$ should be defined globally, which for the
moment means on the whole $A_{00}$, and for this the function $R$ is very useful, because
the computations of Section \ref{S2} show that we can take $a = \frac{|\nabla R|}{|\nabla G|}$.
So we have less choice about $a$ than one could believe: we already defined $G$ near
both circles, with $|\nabla G| = 1$ on both circles, and we want to take $a=1$ near
$\d B_0$ and $\d B_1$, so we need to make sure that $|\nabla R| = |\nabla G|$ 
in both regions, so in particular $|\nabla R| =1$ on both circles. This is what we ensured
with \eqref{3.11}, and now the fact that $|\nabla R| = |\nabla G|$ near the circles follows
from the fact that $G$ is harmonic there, the red curves are known and parameterized
in the right way, so $R$ has to be the usual conjugation of $G$ that we don't even need to compute.
\end{remark}

At this point we have two functions $R$ and $G$, where $R$ is entirely determined
by the red curves and the parameterization from \eqref{3.11}, and we have a little bit of latitude for $G$, 
because we only decided about its values in three precise regions.
This is not too shocking, because we could decide to replace $G$ with a composed function like 
$f \circ G$ (with a nice $f$); we can lift the ambiguity by parameterizing one of the red curves 
$\Gamma_\theta$ by $s \in (\frac14,1)$ and deciding
that $G(z) = s$ along the green curve that contains $\Gamma(s)$.

As was just explained, we want to take 
\begin{equation} \label{3.12}
a(z) = |\nabla G(z)|/ |\nabla R(z)| \text{ on } A_{00},
\end{equation}
because this is our way of making sure that
\begin{equation} \label{3.14}
\dv a \nabla G = 0 \text{ on } A_{00},
\end{equation}
even in the strong sense (because $G$ is smooth).
Let us check that 
\begin{equation} \label{3.14a}
C^{-1} \leq a(z) \leq C \ \text{ on } A_{00}.
\end{equation}
Away from $0$, this is clear because all our functions are smooth,
and we can ensure lower bounds on their gradients. In the small ball $B = B(0,\rho_2)$ 
near $0$, we have a precise formula for $G$, which determines the red curves
but not how they are labelled. This is what gives the picture of Figure \ref{f4},
copied from a piece of Figure \ref{f1}. There is a good choice $R_0$ for $R$, coming from 
the formula \eqref{2.6} and that yield $a=1$, but multiplying $R_0$ by a constant
will only lead to a different constant value of $a$, which is fine too. In fact, if we do not
pay attention, the values of $R$ when we enter $B$ will (by smoothness) be equivalent 
to the values of $R_0$, in the sense that $C^{-1}|\nabla R_0| \leq  |\nabla R| \leq C |\nabla R_0|$.
This will be enough for ellipticity, but we also promised that $a$ would be continuous, so
we modify slightly the way our red curves $\wt \Gamma_\theta$ were organized near 
$I_2$ (or $I_4$) to make sure that for those curves that enter  $B(0,\rho_2)$,
the parameterization speed is proportional to what we would get for $R_0$. 
This way we obtain that 
\begin{equation} \label{3.14b}
a \ \text{ is constant on } A_{00} \cap B(0,\rho_2),
\end{equation}
which is a brutal way to ensure that $a$ will be continuous at $0$. Let us also record that
\begin{equation} \label{3.13}
a(z) = 1 \text{ in }  A_{00} \cap \Big((B_0 \sm B(0,\rho)) \cup (B(z_1, \rho_1) \sm B_1) \Big).
\end{equation}

\begin{figure} [!h] 
\centering
\includegraphics[width=5.cm]{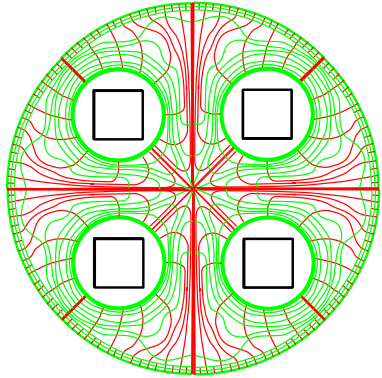}
\caption{The level and gradient lines of $G$ in $A_0$, completed by symmetry}
 \label{f7}
\end{figure}

\ms
This ends our construction of $R$ and $G$ on $A_{00}$. 
We extend $G$ and $a$ to $A_0$ by symmetry. Notice that 
$a$ is continuous on $A_0$; it is possible that we could make it smooth,
by choosing the red curves near $\d A_{00}$ even more carefully so that
the normal derivative of $a$ vanishes there, but we were not courageous enough to try.
Also, $G$ is $C^1$ across the the straight part of the boundary of $\d A_{00}$,
because its normal derivative vanishes (this part of the boundary is a red curve).
Then it is also of class $C^2$, and by symmetry (or reflection principle), \eqref{3.14}
still holds in $A_{0}$. What happens near $0$ is even simpler, because $G$
is harmonic and \eqref{3.14b} holds there.
At this point we have a picture that looks like  Figure \ref{f7}. 

It is now time to glue $G$ to itself in a fractal way.
For each generation $k \geq 0$, denote by $\cQ(k)$ the set of $4^k$
squares that compose $K_k$. For each square $Q \in \cQ(k)$ we denote by
$z(Q)$ the center of $Q$ and $F_Q$ the obvious affine mapping
that sends $Q$ to $Q_0 = K_0$. Thus
\begin{equation} \label{3.15}
F_Q(z) = 4^k (z-z(Q)) \ \text{ for } z \in \R^2.
\end{equation}
We define an annular region $A(Q) = F_Q^{-1}(A_0)$, notice that the 
$A(Q)$, $Q \in \cup_k \cQ(k)$, form a partition of $B_0 \sm K$,
and define functions on $B_0 \sm K$ by
\begin{equation} \label{3.16}
G(z) = 4^{-k} G(F_Q(z)), \, \text{ and } \,
a(z) = a(F_Q(z)) \ \text{ for } z \in A(Q), Q \in \cQ(k).
\end{equation}
Notice that  
\begin{equation} \label{3.17}
\text{$4^{-k-1} \leq G \leq 4^{-k}$ on $A(Q)$ when $Q \in \cQ(k)$,}
\end{equation}
with 
\begin{equation} \label{3.18}
\text{$G(z) = 4^{-k}$ on the exterior boundary $F_Q^{-1}(\d B_0)$,}
\end{equation}
which is also (when $k \geq 1$) the interior boundary of the cube of the previous 
generation that contains $Q$. That is, $G$ is continuous across the boundaries,
and we even claim that
\begin{equation} \label{3.19}
\text{$G$ is harmonic across $F_Q^{-1}(\d B_0)$}.
\end{equation}
This requires a small verification of normal derivatives, or we can directly 
use the formulas. Let us just do the verification when $Q$ is the cube
$Q_1$ of generation $1$; in the general case, we would need to compose both functions
with an additional affine transform. Near $\d B_1$, but at the exterior, we decided in \eqref{3.8}
that $G(z) = \frac14 + \frac14 \ln(4|z-z_1|)$. But on the interior part, we use the formula
\begin{equation} \label{3.20}
G(z) = \frac14 G(F_{Q_1}(z)) = \frac14 G(4(z-z_1)) = \frac14 (1 + \ln(4|z-z_1|)).
\end{equation}
by \eqref{3.7}. This is the same formula (we made it on purpose!), so $G(z)$
is indeed harmonic across the circle.

\begin{figure} [!h]  
\centering
\includegraphics[width=5.cm]{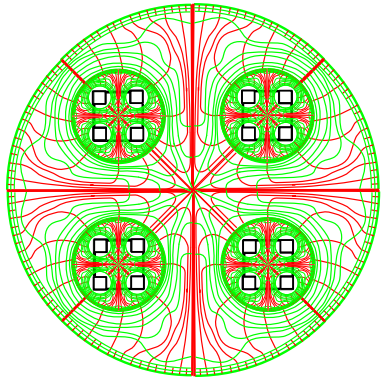}
\caption{The level and gradient lines of $G$ on a larger region $A_0$, completed 
by fractality}
\label{f8}
\end{figure}

We finally obtained our scale invariant elliptic coefficient $a$, and a function 
$G$ which is $L$-harmonic (i.e., satisfies the equation $\dv a \nabla u = 0$)
on $B_0 \sm K$, and is equivalent to $\dist(z,K)$ (by \eqref{3.17}). 
This is essentially all we needed. For the fun of it, we check on Figure \ref{f8} 
that when we put four copies of Figure \ref{f8}, 
reduced by a factor of $4$, in their
correct place in the $B_j$, we get a description of $G$ in a larger region that looks coherent.
With our additional patching law, the continuations of the red curves arrive to $16$ small circles
with the same distances as when they left $\d B_0$, even though the picture seems to say
something different.

In the next section we finally state the main theorem and explain why we essentially 
proved it already. 

\section{A statement of the main theorem}
\label{S4}

Let us recall some of the notation for the main theorem.
Let $K \subset \R^2$ be the Garnett-Ivanov Cantor set 
that was described near \eqref{3.1}, and
let $\mu$ denote the canonical probability measure on $K$, defined by the fact
that $\mu(Q) = 4^{-k}$ for each of the squares of sidelength $4^{-k}$
that compose $K_k$ (so that we don't even need to define $\H^1$).

For each operator $L =  \dv a \nabla$, with $a$ measurable and such that
\begin{equation} \label{4.1}
C^{-1} \leq a(z) \leq C  \, \text{ for } z \in \R^2 \sm K
\end{equation}
(so that $L$ is elliptic), and each point $z \in \R^2 \sm K$, there is a probability measure
$\omega_L^z$ on $K$, called the elliptic measure for $L$ on $K$, and which for instance allows 
one to solve the Dirichlet problem. We can also let $z$ tend to $\infty$, and by an argument that
uses the comparison principle, $\omega_L^z$ 
tends to a probability measure $\omega_L^{\infty}$, which we call the harmonic measure 
with pole at $\infty$. We introduce it here because it gives a cleaner statement.

\begin{theorem} \label{t41}
Retain the notation above.
There exists an absolute constant $C \geq 1$
and a continuous function $a : \R^2 \sm K \to (0,+\infty)$, 
such that \eqref{4.1} holds,  $a(z) = 1$ for $z \in \R^2 \sm B(0,1)$,
and if $\omega_L^\infty$ denotes the harmonic measure, with pole at $\infty$,
associated to the operator $L = \dv a \nabla$ on the domain $\R^2 \sm K$, 
then $\omega_L^\infty = \mu$.
Also,
\begin{equation} \label{4.3}
C^{-1} \mu(A) \leq \omega_L^z(A) \leq C  \mu(A)
\text{ for $z \in \R^2 \sm B(0,1)$ and $A \subset K$ measurable.}
\end{equation}
\end{theorem}

We shall discuss a few variants and improvements in the next section. In the 
meantime we prove the theorem. Let $a$ be the function that was constructed in the last section,
which we extend by taking $a = 1$ on $\R^2 \sm B_0$. Then let $G$ be 
the $L$-harmonic function above, which we extend by taking 
$G(z) = 1 + \ln(|z|)$ for $|z| > 1$. That is, we just extend the formula \eqref{3.7},
and of course $G$ is harmonic on $\R^2 \sm B_0$. We claim that $G$ is 
(a constant times) the Green function for $L$, with the pole at $\infty$, essentially because it is
$L$-harmonic, positive, and vanishes at the boundary.
Here it is even equivalent to $\dist(z,K)$ near $K$, 
but the usual H\"older continuity would be enough.

Now we don't want to use the usual relation between the normal derivative of the Green 
function and the Poisson kernel directly on $K$, because $K$ is irregular, so the simplest seems 
to approximate $K$ by the circular variant of $K_k$ and compute there.

Denote by $D_k = \cup_{Q \in \cQ(k)} F_Q^{-1}(\ol B_0)$ this approximation. Here we use
the notation near \eqref{3.15}, $D_k$ is just a union of balls of radius $4^{-k}$ with the
same centers as the pieces of $K_k$, and $E_k = \d D_k$ is the corresponding
union of circles. Notice that $K \subset D_k$,
$\Omega_k = \R^2 \sm D_k$ is a nice smooth domain, 
$G$ is still $L$-harmonic on $\Omega_k$, and $G_k = G - 4^{-k}$  is clearly 
(a multiple of) the Green function on $\Omega_k$, with pole at $\infty$. 
We can even compute the normal derivative $g_k = \frac{\d G_k}{\d n}$ on 
$\d D_k$, because we have the explicit formulas \eqref{3.16} and \eqref{3.8}.
We find that $g_k$ is a constant, that does not even depend on $k$, so
the ellipitic measure at $\infty$ associated to $L$ on $\Omega_k$
(call it $\omega_k$)  is equal to the invariant measure $\mu_k$ on $\d D_k$. 
We write the reproducing formula
$f(\infty) = \int_{E_k} f d\omega_k = \int_{E_k} \mu_k$ 
for continuous bounded $L$-harmonic functions, 
let $k$ tend to $+\infty$, observe that $\mu_k$ tends to $\mu$ 
(look at the effect on any finite union of sets $F_Q^{-1}(\ol B_0)$),
prove that $\mu$ gives an good reproducing formula on $K$,
and conclude.

The second part of our statement is an easy consequence of the first one
and the change of poles formula (or the comparison principle); if the reader does not like
to take a pole at $\infty$, they can also look at the argument above, notice that
our function $G$ is equivalent to the Green function with the pole $X$, and then derive
\eqref{4.3} with the same argument as sketched above for the Green function at infinity.
\qed

Here we decided to go for the simplest statement, but variants and extensions are possible,
some of which we describe in the next section.

\section{Variants and extensions}
\label{S5}

\begin{figure} [!h] 
\centering
\includegraphics[width=4.cm]{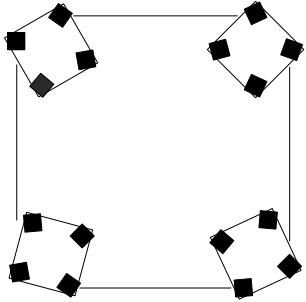}
\caption{The third iteration of a rotating version of the Cantor set}
 \label{f9}
\end{figure}

We start with \ub{rotated Cantor sets}. Suppose that,
in the description of the construction of  $K_k$ from $K_{k-1}$, when we replace each cube $Q$ 
with four squares at its corners, we also allow to rotate $Q$ (or now the four cubes) by any angle 
that depends on $Q$. See Figure \ref{f9} 
for a hint; we leave the reader draw themselves the analogue of Figure \ref{f8} in this case. 
We get a new, apparently more complicated set $K^\ast$; this manipulation makes it much more 
difficult to control the usual harmonic measure on $K^\ast$ (and prove that it is very 
singular), or control the measure of the projections, but here the effect on our construction is just null: 
the domains $A_0$ and $A_{00}$ are different, but the important information is that
we have a (rotation invariant) formula for $G$ on the exterior and interior boundaries,
so that we can glue the pieces. Of course the function $a$ depends on the rotations.

Next consider \ub{differently scaling Green functions}. The simplest way to do this is to take
a function of the function $G(z)$ constructed above, such as $G(z)^{\alpha}$.
This gives a solution of $\wt L = \dv \wt a \nabla$, where 
$\wt a(z) = \alpha a(z) G(z)^{\alpha-1} \sim \dist(z,K)^{\alpha-1}$.
The Green function for $\wt L$ is $G(z)^{\alpha} \sim \dist(z,K)^{\alpha}$. 
This is not too surprising: we trade a different scaling for $G$ against a different scaling
for the degenerate elliptic operator $\wt L$. Notice that in this case too,
the corresponding harmonic measure at $\infty$ is still the invariant measure $\mu$
(what else?).

\ms
There is no special difficulty with replacing $K$ by \ub{different Cantor sets} $K_\beta$
constructed in the same way, but a dimension $\beta \neq 1$. That is, instead
of taking squares of sidelength $4^{-k}$ at the $k$th generation, we take squares 
of sidelength $\lambda^{k}$ for some $\lambda \neq \frac14$. The dimension $\beta$ is limited
by the fact that we want the balls $B_j$, $1 \leq j \leq 4$ to contain the corresponding
cubes $Q_j$ and be disjoint. 
 
 Something interesting happens in the construction, which is related to Remark \ref{r31}.
When we construct the function $G$, we decide that it is equal to $1$
on $\d B_0$, and to some constant $\gamma < 1$ on the analogue of $\d B_1$.
Then we choose the function $R$ so that the transfer map from
$\d B_0$ to $\d B_1$ goes at a constant speed $\tau$. The length to cover is
$\pi/4$ on $\d B_0$ and $\lambda\pi$ on $\d B_1$, so the gradient of $R$
is $4\lambda$ times smaller on $\d B_1$. 

At the same time, if we want to use a fractal formula like \eqref{3.16}, 
we should write it down as
\begin{equation} \label{5.1}
G(z) = \gamma^k G(F_Q(z)), 
\end{equation}
where the derivative of $F_Q(z)$ is now $\lambda^{-k}$. Hence
$|\nabla G|$ is $\gamma \lambda^{-1}$ times larger on $\d B_1$ than 
on $\d B_0$ and $a = |\nabla R|/|\nabla G|$ is $4\gamma$ times larger on on $\d B_1$.

If we want to construct an elliptic operator, we should be able to get $a=1$
on both circles, so $\gamma = 1/4$. This means
that for points of generation $k$, the distance to $K$ is like
$\lambda^k$, but the size of the Green function is now like $\gamma^k = 4^{-k}$,
which is also roughly the measure in $K$ of a ball of radius $\lambda^k$.
In other words, $G(z) \sim \dist(z,K)^{d}$, where $d$ is the dimension of the Cantor set.

In general, we can pick $\gamma \neq 1/4$, get a Green function with essentially
any other homogeneity, but then our operator $L = \dv A \nabla$ will be degenerate elliptic,
with $|a(z)| \sim \dist(z,K)^{\alpha}$ for some $\alpha \neq 0$.

We can pursue all this a little further, and replace the balls in the construction of $G$
with objects with a different shape (in fact our first constructions were like that)
but then we are no longer allowed to rotate the squares.

\begin{figure} [!h]
\centering
\includegraphics[width=4.cm]{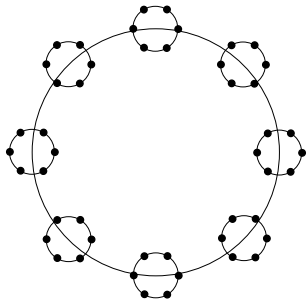}
\caption{The third iteration of a variable scale/multiplicity analogue of $K$;
there is no point in trying to draw polygones in this case}
\label{f10}  
\end{figure}

In fact (as in \cite{Bat}, for instance), we can let the dilation ratio depend on the scale
(as long as we keep some uniformity), or take Cantor sets based on dividing each
box into more that $4$ pieces, and probably combine all of the above.
See Figure \ref{f10} for a hint. 

It would be nice to have an operator like the one above, with the special form
\eqref{1.1} (otherwise this is too easy), associated to the standard snowflake,
or some general Reifenberg flat domains. This is tempting, but one would have to be careful with
the construction and the verification that $C^{-1} \leq a \leq C$.
We can also ask the same questions in higher dimensions, i.e., for instance for a product of three
Cantor sets of dimension $1/3$ in $\R^3$.

\end{document}